\documentclass[12pt]{article} 
\usepackage{amsmath,amsthm, amssymb} 
\usepackage{amssymb,latexsym}

\newtheorem*{theorem}{Theorem}

\newtheorem*{cor}{Corollary}

\begin{document}

\baselineskip=17pt 

\title{\bf An identity for the Kloosterman sum}

\author{\bf D. I. Tolev}

\date{}
\maketitle

\begin{abstract}
We establish a simple identity and using it we find a new proof of a result of Kloosterman.

\smallskip

\end{abstract}

The Kloosterman sum is the defined by
\begin{equation} \label{5}
K(p; a, b) = \sum_{x=1}^{p-1} e_p \left( a x + b \overline{x} \right) ,
\end{equation}
where $p$ is a prime, $a$ and $b$ are integers,
$\overline{x}$ is the inverse of $x$ modulo $p$
and $e_p(\alpha) = \exp \left( \frac{2 \pi i \alpha}{p} \right)$.
It is clear that it takes always real values. This sum was introduced in 1926 by Kloosterman~\cite{Kloos}  
and he established that
\begin{equation} \label{10}
 |K(p; a, b)| \le 3^{1/4} \, p^{3/4} \qquad \text{for} \qquad p \nmid ab .
\end{equation}

In 1948 A.Weil~\cite{Weil} improved substantially the estimate \eqref{10} and obtained the following deep and important 
inequality:
\begin{equation} \label{20}
 |K(p; a, b)| \le 2 \, \sqrt{p} \qquad \text{for} \qquad p \nmid ab .
\end{equation}
Later  Stepanov~\cite{Step} found an elementary proof of \eqref{20}
(see also Iwaniec and Kovalski~\cite{IwKo}, Chapter~11).
Information about the applications of Kloosterman's sum in analytic number theory as well as a simple proof of 
\eqref{10} can be found in Heath-Brown's paper~\cite{HB}. 
Another proof of \eqref{10} is available in the recent preprint 
\cite{FGK} from Fleming, Garcia and Karaali.

\bigskip

In this short note we present an identity for the Kloosterman sum and using it we find
a new proof of \eqref{10} (with smaller constant in the right-hand side of this inequality). 

\bigskip

From this point onwards we assume that $p>2$ is a fixed prime 
and let $\left( \frac{\cdot}{p} \right)$ be the Legendre symbol. We write for simplicity
$K(a, b) = K(p; a, b)$. Our result is the following
\begin{theorem}
For any integers $a$, $b$ such that $p \nmid ab$ we have
\begin{equation} \label{25}
    K(a, b)^2 = p + \sum_{l=1}^p \left( \frac{l^2 - 4 l}{p} \right) \, 
    K(a, lb) .
\end{equation}
\end{theorem}

{\bf Proof:}
Using \eqref{5} we find
\begin{align}
  K(a, b)^2 
    & = \sum_{1 \le x, y \le p-1} e_p \left( a \, (x - y) + b \, (\overline{x} - \overline{y}) \right)
      = 
   \sum_{h=1}^p  e_p (ah) \sum_{\substack{1 \le x, y \le p-1 \\ x-y \equiv h \pmod{p}}}
  e_p \left(  b \, (\overline{x} - \overline{y}) \right)  \notag \\
    & = p-1 + Y(a, b) , 
    \label{30}
\end{align}
where
\[
  Y(a, b) =  \sum_{h=1}^{p-1}  e_p (ah) 
    \sum_{\substack{1 \le y \le p-1 \\ p \, \nmid \, y+h }}
  e_p \left(  b \, \left( \overline{y+h} - \overline{y} \right) \right) .
\]
We put $y = hz$ in the inner sum and obtain
\[
  Y(a, b) =  \sum_{h=1}^{p-1}  e_p (ah) 
    \sum_{z=1}^{p-2}
  e_p \left(  b \, \overline{h} \, \left( \overline{z+1} - \overline{z} \right) \right) .
\]
Now we change the order of summation and use \eqref{5} to get
\begin{equation} \label{50}
 Y(a, b) = \sum_{z=1}^{p-2} K \left(a, b \left( \overline{z+1} - \overline{z} \right) \right) 
   = \sum_{l=1}^{p-1} K(a, lb) \; \lambda_l ,
\end{equation}
where $\lambda_l$ is the number of integers $z$ such that
$ 1 \le z \le p-2 $ and $\overline{z+1} - \overline{z} \equiv l \pmod{p} $.
We easily  see 
that $\lambda_l$ equals the number of solutions of the congruence
$lz^2 + lz + 1 \equiv 0 \pmod{p}$, hence from the properties of the Legendre symbol it follows that
\begin{equation} \label{60}
  \lambda_l = 1 + \left( \frac{l^2 - 4l}{p} \right) .
\end{equation}
From \eqref{5} and our assumption $p \nmid ab$ we get
\begin{equation} \label{65}
 \sum_{l=1}^{p-1} K (a, lb) = \sum_{x=1}^{p-1} e_p (ax) 
     \sum_{l=1}^{p-1} e_p \left( bl \overline{x}\right) = - \sum_{x=1}^{p-1} e_p (ax) = 1 .
\end{equation}
The identity \eqref{25} is a consequence of \eqref{30} -- \eqref{65}.

\hfill
$\square$
 
\bigskip

Now we obtain immediately the following 
\begin{cor}
If $p \nmid ab$ then
\begin{equation} \label{70}
  |K(a, b)| \le \sqrt{ p + p^{3/2} } .
\end{equation}
\end{cor}

{\bf Proof:}
Denote by $Z$ the second term in the right-hand side of \eqref{25}. From Cauchy's inequality 
we get
\[
  |Z| \le  p^{1/2} \, \left( \sum_{l=1}^p K^2(a, lb) \right)^{1/2} = p^{1/2} \, Z_1^{1/2} ,
\]
say. From \eqref{5} it follows that
\[
  Z_1   = \sum_{1 \le x, y \le p-1}  e_p \left( a (x-y) \right) 
     \sum_{l=1}^p e_p \left( bl \left( \overline{x} - \overline{y} \right) \right) = p \, (p-1) .
\]
Hence $|Z| \le p^{3/2}$ and using \eqref{25} we obtain \eqref{70}.

\hfill
$\square$

\bigskip

Finally we mention that by the same method we can estimate also the sum
\[
   K_r(p; a, b) = \sum_{x=1}^{p-1} e_p \left( a x^r + b \overline{x} \right) ,
\]
for arbitrary positive integer $r$. We can prove that
$  K_r(p; a, b) \ll_r p^{3/4} $ for $p \nmid ab$,
but we shall not give the details here.

\bigskip

{\bf Acknowledgments:} 
The present research was supported by Sofia University Grant 172/2010.

\bigskip
\bigskip

\vbox{
\hbox{Faculty of Mathematics and Informatics}
\hbox{Sofia University ``St. Kl. Ohridsky''}
\hbox{5 J.Bourchier, 1164 Sofia, Bulgaria}
\hbox{ }
\hbox{Email: dtolev@fmi.uni-sofia.bg}}

\end{document}